\documentclass{article}
\usepackage{fullpage}
\usepackage{amsmath,amssymb,eqnarray}
\usepackage{fancybox,graphicx,epsf,color,array}
\usepackage{enumitem}   

\renewcommand{\b}[1]{{\bf #1}}

\newcommand{\fz}[3]{#1:\, #2 \rightarrow #3}

\renewcommand{\r}[1]{(\ref{#1})}
\newcommand{\bi}{\begin{itemize}}
\newcommand{\ei}{\end{itemize}}
\newcommand{\be}{\begin{enumerate}}
\newcommand{\ee}{\end{enumerate}}

\newcommand{\bd}{\begin{description}}
\newcommand{\ed}{\end{description}}
\renewcommand{\i}{\item}

\newcommand{\bqn}{\begin{eqnarray}}
\newcommand{\eqn}{\end{eqnarray}}
\newcommand{\eqnn}{\nonumber\end{eqnarray}}
\newcommand{\eqnl}[1]{\label{#1}\end{eqnarray}}

\newcommand{\nn}{\nonumber\\}

\newcommand{\ba}[1]{\begin{array}{#1}}
\newcommand{\ea}{\end{array}}

\newcommand{\R}{\mathbb{R}}

\newcommand{\N}{\mathbb{N}}
\newcommand{\Z}{\mathbb{Z}}
\newcommand{\Q}{\mathbb{Q}}

\newcommand{\bproof}{\begin{proof}}
\newcommand{\eproof}{\end{proof}}
\newtheorem{Theorem}{\bf Theorem}
\newtheorem{lemma}[Theorem]{\bf Lemma}
\newtheorem{corollary}[Theorem]{\bf Corollary}
\newtheorem{definition}[Theorem]{\bf Definition}
\newtheorem{proposition}[Theorem]{\bf Proposition}
\newtheorem{remark}[Theorem]{\bf Remark}
\newenvironment{proof}{{\it Proof.}~~}{\hfill$\square$}

\newcommand{\bt}{\begin{Theorem}}
\newcommand{\et}{\end{Theorem}}
\newcommand{\bl}{\begin{lemma}}
\newcommand{\el}{\end{lemma}}
\newcommand{\bp}{\begin{proposition}}
\newcommand{\ep}{\end{proposition}}
\newcommand{\bc}{\begin{corollary}}
\newcommand{\ec}{\end{corollary}}
\newcommand{\bdeff}{\begin{definition}}
\newcommand{\edeff}{\end{definition}}
\newcommand{\brem}{\begin{remark}\rm}
\newcommand{\erem}{\end{remark}}

\newcommand{\eps}{\varepsilon}
\newcommand{\ga}{\gamma}

\newcommand{\Pt}[1]{\left( #1 \right)}
\newcommand{\Pg}[1]{\left\{ #1 \right\}}
\newcommand{\Pq}[1]{\left[ #1 \right] }

\newcommand{\Pabs}[1]{\left| #1 \right|}

\newcommand{\schema}[1]{\b{\sc #1}}

\newcommand{\M}{\mathcal{M}}
\renewcommand{\P}{\mathcal{P}}

\newcommand{\W}{\mathcal{W}}

\newcommand{\supp}{\mathrm{supp}}

\newcommand{\Mu}{\mathcal{M}(\R^d)}
\newcommand{\weak}{\rightharpoonup}

\newcommand{\diam}{\mathrm{diam}}

\newcommand{\gw}{W^g}
\renewcommand{\W}{\mathcal{W}}
\newcommand{\GW}{\W^g}
\newcommand{\A}{\mathcal{A}}

\title{Measure dynamics with Probability Vector Fields and sources}

\begin{document}

\maketitle

\centerline{\scshape Benedetto Piccoli}
\medskip
{\footnotesize
 \centerline{Department of Mathematical Sciences, Rutgers University - Camden }
   \centerline{ Camden, NJ, USA.}
   \centerline{{\tt piccoli@camden.rutgers.edu}}
} 

\medskip

\centerline{\scshape Francesco Rossi}
\medskip
{\footnotesize
 \centerline{ Dipartimento di Matematica ``Tullio Levi--Civita''}
   \centerline{Universit\`a degli Studi di Padova}
   \centerline{Padova, Italy}
   \centerline{{\tt francesco.rossi@math.unipd.it}}
}

\bigskip


\begin{abstract} We introduce a new formulation for differential equation describing dynamics of  measures on an Euclidean space, that we call Measure Differential Equations with sources. They mix two different phenomena: on one side, a transport-type term, in which a vector field is replaced by a Probability Vector Field, that is a probability distribution on the tangent bundle;  on the other side, a source term. Such new formulation allows to write in a unified way both classical transport and diffusion with finite speed, together with creation of mass.

The main result of this article shows that, by introducing a suitable Wasserstein-like functional, one can ensure existence of solutions to Measure Differential Equations with sources under Lipschitz conditions. We also prove a uniqueness result under the following additional hypothesis: the measure dynamics needs to be compatible with dynamics of measures that are sums of Dirac masses.
\end{abstract}

\vspace{1cm}

{{\bf Keywords:} Measure dynamics, Probability Vector Fields, Wasserstein distance, generalized Wasserstein distance}

{{\bf MSC2010}: 35S99, 35F20, 35F25} 



\section{Introduction}

The problem of optimal transportation, also called  Monge-Kantorovich problem, has been intensively studied in the mathematical community. Related to such problem, the definition of the Wasserstein distance in the space of probability measure has revealed itself to be a powerful tool, in particular for dealing with dynamics of measures (like the transport PDE, see e.g. \cite{gradient}). For a complete introduction to Wasserstein distances, see \cite{villani,old-new}.

This approach has at least two main limits. The first is that the use of transport equation, together with their counterpart in terms of ordinary differential equations \cite{amb-inv,ACF-arma}, does not allow to model neither mass diffusion nor concentration phenomena. The second one is that the Wasserstein distance $W_p(\mu,\nu)$ is defined only if the two measures $\mu,\nu$ have the same mass, then PDEs with sources cannot be studied with such tools.

Both limits were recently overcome by a variety of contributions. The first was addressed in \cite{MDE}, in which a generalization of the concept of vector fields was introduced. Such tool, called Probability Vector Field (PVF in the following), allows to model concentration and diffusion phenomena in the formalism of the transport equation, then being able to translate several useful techniques from dynamical system.

The second limit was addressed by a series of papers introducing generalizations of the Wasserstein distance to measures with different masses. In \cite{gw} we defined a generalized Wasserstein distance $\gw(\mu,\nu)$, combining the standard Wasserstein and $L^1$ distances. In rough words, for $\gw(\mu,\nu)$ an infinitesimal mass $\delta\mu$ of $\mu$ (or $\nu$) can either be removed at cost $a |\delta\mu|$, or moved from $\mu$ to $\nu$ at cost $bW_p(\delta\mu,\delta\nu)$. This distance is a generalization of the so-called {\it flat distance}. Other generalizations of the Wasserstein distance, with the same spirit of allowing sources of mass, are studied in \cite{chizat,KMV16,liero1,liero2}. As a consequence, sources of mass can be introduced in the transport equation, even when they depend on the measure itself, see \cite{gw2}.\\

The goal of this article is to define a new class of equations,
which are able to describe complex dynamics in the space of measures,
including mass diffusion, concentration and sources.
The idea is to merge two different dynamics, already individually described in \cite{MDE,gw}, and couple them.\\
The first contribution is given by dynamics induced by Probability Vector Fields (PVF in the following), recently introduced in \cite{MDE}. There, the equation
\bqn
\dot \mu=V[\mu]
\eqnl{e-MDE}
is considered, where $V:\P(\R^n)\to \P(T\R^n)$ is a function from the space $\P(\R^n)$ of probability measures to the space $\P(T\R^n)$ of probability measures of the tangent space $T\R^n$. The idea of such function is to describe the infinitesimal spreading of the mass $\mu(x)$ in a point $x$ along the velocities described by the measure $V[\mu](x,\cdot)$ on the fiber $T_x\R^n$. Given the projection $\pi:T\R^n\to\R^n$ defined by $\pi(x,v)=x$, we also require  $\pi\#V[\mu]=\mu$, i.e. that the projection of $V[\mu]$ from $\P(T\R^n)$ to $\P(\R^n)$ coincides with $\mu$. This is the measure counterpart of the fact that a vector field is a section of the tangent bundle. The main contribution of \cite{MDE} is to introduce conditions ensuring existence and/or uniqueness of the solution of the Cauchy problem with dynamics \r{e-MDE}. In particular, two key tools are defined: the first is a new non-negative operator $\W$, based on the Wasserstein distance and enjoying some of its properties, on the space $\P(T\R^n)$. The idea is that $\W$ measures the cost of the minimizing transference plan on fibers, among plans whose projections are optimal on the base space. The formal definition is given in Definition \ref{d-W}. If one assumes that $V$ from $\P(\R^n)$ endowed with the Wasserstein distance to $\P(T\R^n)$ endowed with $\W$ is Lipschitz, then there exists at least one solution to \r{e-MDE}. The second tool is the definition of Dirac germs, that are specific choices of solutions to \r{e-MDE} for measures composed of Dirac deltas only. Fixed a Dirac germ for \r{e-MDE}, then for each initial measure there exists at most one solution to \r{e-MDE} that is compatible with such chosen germ. In some specific but relevant cases, the coupling of Lipschitz continuity of $V$  with the choice of a compatible Dirac germ ensures both existence and uniqueness of a solution to \r{e-MDE}.\\
The second contribution is given by sources and sinks. In this case, the dynamics reads as
\bqn
\dot \mu= s,
\eqnl{e-source}
where $s$ is a measure on the space $\R^n$, representing a source/sink of mass. The description of such Partial Differential Equation with a fixed source $s$ is very classical, since the solution is clearly $\mu_t=\mu_0+ts$. Instead, we introduced in \cite{gw} new conditions to ensure that the dynamics \r{e-source} is well posed even when the source $s[\mu]$ depends on the whole measure $\mu$ itself. The key tool is the introduction of a new distance on the space of measures with finite mass, called the generalized Wasserstein distance $\gw$. If $s$ is Lipschitz with respect to this distance, then one has existence and uniqueness of the solution to the Cauchy problem with dynamics \r{e-source}. 

For simplicity, from now on we restrict ourselves to the space $\M(\R^d)$ of Borel measures with bounded support and finite mass. In this space, the generalized Wasserstein distance $\gw(\mu,\nu)$ is always finite, while the standard Wasserstein distance $W(\mu,\nu)$ is defined only if the masses of the two measures coincide, i.e. $\mu(\R^n)=\nu(\R^n)$. We endow the space $\M(\R^n)$ with the topology of weak convergence; this coincides with the topology induced by the generalized Wasserstein distance, see Proposition \ref{p-topology} below.\\
We are now ready to define Measure Differential Equations with Source:
\bqn
\dot \mu =V[\mu]\oplus s[\mu],
\eqnl{e-MDES}
where $V[\mu]$ is a PVF $V:\M(\R^n)\to \M(T\R^n)$ and $s[\mu]$ is a source $s:\M(\R^n)\to \M(\R^n)$. The goal is to prove existence and/or uniqueness of a solution to the associated Cauchy problem, under the joint hypotheses ensuring existence and/or uniqueness for each of the dynamics \r{e-MDE} and \r{e-source}. 
More precisely, we first give the definition of a solution to \r{e-MDES}:
\bdeff[Solution to \r{e-MDES}] \label{d-sol} A solution to \r{e-MDES} is a continuous curve $\mu:[0,T]\to \M(\R^n)$ satisfying the following condition: for each $f\in C^\infty_c(\R^n)$
\bi
\i the integral $\int_{T\R^n} (\nabla f(x)\cdot v)\, dV[\mu(\tau)](x,v)$ is defined for almost every $\tau\in[0,T]$;
\i the map $\tau\to \int_{\R^n} f(x) ds[\mu(\tau)](x)$ belongs to $L^1([0,T])$;
\i the map $t \to \int f\, d\mu(t)$ is absolutely continuous, and it satisfies 
\bqn\frac{d}{dt} \int_{\R^n} f\,d \mu(t)=\int_{T\R^n} (\nabla f(x)\cdot v)\, dV[\mu(t)](x,v)+\int_{\R^n} f(x)ds[\mu(t)](x)
\eqnl{e-weak}
for almost every $t\in[0,T]$.
\ei
\edeff
Such definition is pretty weak and can not allow uniqueness results, thus we are also interested in stronger properties for solutions to \r{e-MDES}. In particular, we focus on existence of semigroups of solutions, whose definition in this setting is given below.
\bdeff \label{d-semigroup} A Lipschitz semigroup $S_t$ of solutions to \r{e-MDES} is a map $S:[0,T]\times \M(\R^n)\to \M(\R^n)$ satisfying:
\be
\i $S_0\mu=\mu$ and $S_{t+s}\mu=S_tS_s \mu$;
\i the map $t\to S_t\mu$ is a solution to \r{e-MDES};
\i for every $R,M>0$ there exists $C=C(R,M)>0$ such that if $\supp(\mu)\cup\supp(\nu)\subset B(0,R)$ and $\mu(\R^n)+\nu(\R^n)\leq M$, it then holds
\be
\i $\supp(S_t\mu)\subset B(0,e^{Ct}(R+M+1))$;
\i $\gw(S_t\mu,S_t\nu)\leq e^{Ct}\gw(\mu,\nu)$;
\i $\gw(S_t\mu,S_s\mu)\leq C |t-s|$.
\ee
\ee
\edeff
We also need to define a natural tool, merging properties of the operator $\W$ on $\P(T\R^n)$ with the setting of the generalized Wasserstein distance $\gw$ on $\M(\R^n)$. Such non-negative operator, that we denote by $\GW$, measures the minimal standard Wasserstein distance on the fiber between transference plans whose projections give a minimizing decomposition for the generalized Wasserstein distance on the base space. The operator is precisely defined in Section \ref{s-GW}.

We are now ready to state the two main results of this article. The first deals with existence of a solution to \r{e-MDES}, while the second focuses on uniqueness.
\bt \label{t-existence} Consider the Measure Differential Equation with Source \r{e-MDES} with the following two sets of hypotheses:
\bd
\i[(V)] The Probability Vector Field $V:\M(\R^n)\to \M(T\R^n)$ satisfies:
	\bd
	\i[(V1)] support sublinearity: there exists $C>0$ such that for all $\mu\in \M(\R^n)$ it holds
	$$\sup_{(x,v)\in \supp(V[\mu])} |v|\leq C(1+\sup_{x\in \supp(\mu)} |x|);$$
	\i[(V2)] Lipschitz continuity: for each $R>0$ there exists $K=K(R)>0$ such that $\supp(\mu)\cup\supp(\nu)\subset B(0,R)$ implies
	\bqn
	\GW(V[\mu],V[\nu])\leq K \gw(\mu, \nu);
	\eqnl{e-VLip}
	\ed
\i[(s)] The source $s:\M(\R^n)\to \M(\R^n)$ satisfies:
	\bd
	\i[(s1)] Lipschitz continuity: there exists $L$ such that for all $\mu,\nu\in \M(\R^n)$ it holds
\bqn\gw(s[\mu],s[\nu])\leq L \gw(\mu,\nu);\eqnl{e-sLip}
	\i[(s2)] uniform boundedness of the support: there exists $R$ such that for all $\mu\in \M(\R^n)$ it holds $\supp(s[\mu])\subset B_R(0)$.
	\ed
\ed
Then, there exists a Lipschitz semigroup of solutions to \r{e-MDES} in the sense of Definition \ref{d-semigroup}.
\et
\bt \label{t-uniqueness} Consider the Measure Differential Equation with Source \r{e-MDES} satisfying Hypotheses {\bf (V1), (s)} recalled in Theorem \ref{t-existence}. Choose a Dirac germ $\gamma$, as in Definition \ref{d-germ} below. Then, there exists at most one Lipschitz semigroup compatible with $\gamma$, in the sense of Definition \ref{d-compat2} below.
\et

Several corollaries about existence and/or uniqueness of the solutions to \r{e-MDES} can be directly derived from corresponding results about PVFs from \cite{MDE}. In particular, one can observe that the uniqueness property depends on the PVF $V$ only, and not on the source $s$. We then have the two following remarkable cases:
\bi
\i Let $V[\mu]=\mu\otimes \delta_{v(x)}$ with $v$ locally Lipschitz vector field with sub-linear growth. Then, \r{e-MDES} admits a unique Lipschitz semigroup, obtained as the limit of the discretization described in Section \ref{s-existence}.
\i Fix $\phi:[0,+\infty)\to\R$ an increasing function. In the space $\R$, define
$V_\phi[\mu]=\mu\otimes J_\phi(x)$, where
\bqn
J_\phi(x)=\begin{cases}
\delta_{\phi(F_\mu(x))}&\mbox{~~if~}F_\mu(x^-)=F_\mu(x),\\
\frac{\phi\#\Pt{\chi_{[F_\mu(x^-),F_\mu(x)]}\lambda}}{F_\mu(x)-F_\mu(x^-)}&\mbox{~~otherwise,}
\end{cases}
\eqnn
$F_\mu(x)=\mu((-\infty,x])$ is the cumulative distribution of $\mu$, and $\lambda$ is the Lebesgue measure. This choice of the PVF allows to have solutions that diffuse with finite velocities, see \cite[Section 7.1]{MDE} for more details. In this case, for any choice of the source $s$ satisfying {\bf (s)}, one has existence of a solution to \r{e-MDES}. Even though this solution is not unique, in general, there exists a unique semigroup obtained by the limit of the  discretization algorithm described in Section \ref{s-existence}.
\ei

\brem Observe that hypotheses in Theorem \ref{t-existence} are not sharp, in general. For example, in  {\bf (V2)}, the Lipschitz constant $K$ in \r{e-VLip} can depend on $|\mu|$, with the only requirement of having $\sup_{m\in[0,M]}K(m)<+\infty$ for all finite $M$.

Similarly, condition {\bf (s2)} can be replaced by any condition ensuring uniform boundedness of the supports, such as the existence of a radius $R$ such that $\supp(\mu)\subseteq B(0,R')$ with $R'>R$ implies $\supp(s[\mu])\subseteq B(0,R')$.
\erem

The structure of the article is the following. In Section \ref{s-prel} we fix the notation and recall main properties of the tools used later: the Wasserstein distance, the generalized Wasserstein distance and Measure Differential Equations with Probability Vector Fields. In the main Section \ref{s-proof}, we prove the results of this paper. In Section \ref{s-existence}, we prove Theorem \ref{t-existence} about existence of a solution to \r{e-MDES}, while in Section \ref{s-uniqueness}, we prove Theorem \ref{t-uniqueness} about uniqueness.

\section{Dynamics in generalized Wasserstein Spaces} \label{s-prel}

In this section, we fix the notation and define the main tools used in the rest of the article: the Wasserstein distance, the generalized Wasserstein distance and Measure Differential Equations with Probability Vector Fields.

\subsection{The Wasserstein distance}
We use $\Mu$ to denote the space of positive Borel regular measures with bounded support and finite mass on $\R^d$. Given $\mu,\mu_1$ Radon measures (i.e. positive Borel measures with locally finite mass), we write $\mu_1\ll \mu$
if $\mu_1$ is absolutely continuous with respect to $\mu$, while we write $\mu_1\leq\mu$
if $\mu_1(A)\leq\mu(A)$ for every Borel set $A$. We denote by $|\mu|:=\mu(\R^d)$ the norm 
of $\mu$ (also called its mass). More generally, if $\mu=\mu^+-\mu^-$ is a signed Borel measure, we define $|\mu|:=|\mu^+|+|\mu^-|$. 

Given a Borel map $\fz{\gamma}{\R^d}{\R^d}$, the push forward of a measure $\mu\in\Mu$
is defined by:
\bqn
\gamma\#\mu(A):=\mu(\gamma^{-1}(A)).
\eqnn
Note that the mass of $\mu$ is identical to the mass of $\gamma\#\mu$.
Therefore, given two measures $\mu,\nu$ with the same mass, one may look for $\gamma$
such that  $\nu=\gamma\#\mu$ and it minimizes the cost
$$I\Pq{\gamma}:=|\mu|^{-1}\,\int |x-\ga(x)|^p \,d\mu(x).$$
This means that each infinitesimal mass $\delta\mu$ is sent to $\delta \nu$ and that its infinitesimal cost is  the $p$-th power of the distance between them. 
Such minimization problem is known as the Monge problem. A generalization of the Monge problem is defined as follows.
Given a probability measure $\pi$ on $\R^d\times\R^d$, one can interpret $\pi$ as a method to transfer a measure $\mu$ on $\R^d$ to another measure on $\R^d$ as follows: each infinitesimal mass on a location $x$ is sent to a location $y$ with a probability given by $\pi(x,y)$. Formally, $\mu$ is sent to $\nu$ if the following properties hold:
\bqn
|\mu|\,\int_{\R^d} d\pi(x,\cdot)=d\mu(x),\qquad \qquad |\nu|\,\int_{\R^d} d\pi(\cdot,y)=d\nu(y).
\eqnl{e-pi}
Such $\pi$ is called a transference plan from $\mu$ to $\nu$ and we denote the set of transference plans from $\mu$ to $\nu$ by $P(\mu,\nu)$. 
A condition equivalent to \r{e-pi} is that, for all $f,g\in C^\infty_c(\R^d)$ it holds $|\mu|\,\int_{\R^d\times \R^d} (f(x)+g(y))\,d\pi(x,y) = \int_{\R^d} f(x)\,d\mu(x)+ \int_{\R^d} g(y)\,d\nu(y)$.

One can define a cost for $\pi$ as follows $$J\Pq{\pi}:=\int_{\R^d\times\R^d} |x-y|^p \,d\pi(x,y)$$ and look for a minimizer of $J$ in $P(\mu,\nu)$. Such problem is called the Monge-Kantorovich problem. It is important to observe that such problem is a generalization of the Monge problem. The main advantage of this approach is that a minimizer of $J$ in $P(\mu,\nu)$ always exists. We then denote by $P^{opt}(\mu,\nu)$ the set of transference plans that are minimizers of $J$, that is always non-empty.

One can thus define on $\mathcal{M}_1(\R^d)$ the following operator between measures of the same mass, called the \b{Wasserstein distance}:
\bqn
W_p(\mu,\nu)=|\mu|(\min_{\pi\in P(\mu,\nu)} J\Pq{\pi})^{1/p}.
\eqnn
It is indeed a distance on the subspace of measures in $\Mu$ with a given mass, see \cite{villani}. It is easy to prove that $W_p(k\mu,k\nu)= W_p(\mu,\nu)$ for $k\geq 0$, by observing that $ P(k\mu,k\nu)= P(\mu,\nu)$ and that $J\Pq{\pi}$ does not depend on the mass.

From now on, we only consider the Wasserstein distance with parameter $p=1$, that will then be denoted by $W(\mu,\nu)$. It satisfies the following fundamental dual property.

\bp \label{p-dual}[Kantorovich-Rubinstein duality] Let $\mu,\nu\in$$\mathcal{M}_1(\R^d)$. It then holds
\bqn
W(\mu,\nu)=\sup\Pg{\int f d(\mu-\nu)\,\mbox{~~s.t.~~}\  Lip(f)\leq 1}\label{e-dualdef}
\eqn
\ep

Such property plays a crucial role in the theory of PVF, see \cite{MDE}. It is then unclear if a corresponding theory can be generalized to any $p>1$.

\subsection{The generalized Wasserstein distance}

In this section, we provide a definition of the generalized Wasserstein distance, introduced in \cite{gw,gw2}, together with some useful properties. We consider here the generalized Wasserstein distance with parameters $a=b=1$, to simplify the notation, and $p=1$.

\bdeff \label{d-gw}
Let $\mu,\nu\in\Mu$ be two measures. We define the functional
\bqn
\gw(\mu,\nu):=\inf_{\tilde\mu,\tilde\nu\in\Mu,\,|\tilde\mu|=|\tilde\nu|}|\mu-\tilde\mu|+|\nu-\tilde\nu|+W(\tilde\mu,\tilde\nu).
\eqnl{e-gw}
\edeff

We now provide some properties of $\gw$. Proofs can be adapted from those given in \cite{gw}. 
\bp \label{p-base}
The following properties hold:\\
1. The infimum in \r{e-gw} coincides with
$$\inf_{\tilde\mu\leq\mu,\tilde\nu\leq\nu,\,|\tilde\mu|=|\tilde\nu|}|\mu-\tilde\mu|+|\nu-\tilde\nu|+W(\tilde\mu,\tilde\nu),$$
where we have added the constraint $\tilde\mu\leq\mu,\tilde\nu\leq\nu$.\\
2. The infimum in \r{e-gw} is attained by some $\tilde\mu,\tilde\nu$.\\
3. The functional $\gw$ is a distance on $\Mu$.\\
4. It holds $\gw(\mu,0)\leq |\mu|$.\\
5. It holds
\bqn
||\mu|-|\nu||\leq \gw(\mu,\nu)
\label{e-gwL1}.
\eqn
6. If $|\mu|=|\nu|$, it holds
\bqn
\gw(\mu,\nu)\leq W(\mu,\nu)\label{e-gww}
\eqn
\ep

We recall now some useful topological results related to the metric space $\Mu$ when endowed with the generalized Wasserstein distance. We first recall the definition of tightness in this context.
\bdeff \label{d-tight} A set of measures $M$ is tight if for each $\eps>0$ there exists a compact $K_\eps$ such that $\mu(\R^d\setminus K_\eps)<\eps$ for all $\mu\in M$.
\edeff

We now recall the definition of weak convergence of measures, as well as an important result about convergence with respect to the generalized Wasserstein distance, see \cite[Theorem 13]{gw}.

\bdeff Let $\Pg{\mu_n}$ be a sequence of measures in $\R^d$, and $\mu$ a measure. We say that $\mu^n$ converges to $\mu$ with respect to the weak topology, and we write $\mu_n\weak \mu$, if for all functions $f\in C^\infty_c$ it holds
$$\lim_{n\to\infty} \int f\,d\mu_n=\int f\,d\mu.$$
\edeff

\bp \label{p-topology}
Let $\Pg{\mu_n}$ be a sequence of measures in $\R^d$, and $\mu\in\Mu$. Then
$$\gw(\mu_n,\mu) \to 0\mbox{~~~~~~is equivalent to~~~~~~}\mu_n\weak \mu \mbox{~~and~~}\Pg{\mu_n}\,\mbox{is tight}.$$
\ep

We also recall the result of completeness, see \cite[Proposition 15]{gw}.
\bp \label{p-complete}
The space $\Mu$ endowed with the distance $\gw$ is a complete metric space.
\ep

The generalized Wasserstein distance also satisfies a useful dual formula, showing that it coincides with the so-called flat distance. See \cite{gw2}.
\bp \label{p-dualgw} Let $\mu,\nu\in\Mu$. It then holds
\bqn
\gw(\mu,\nu)=\sup\Pg{\int f d(\mu-\nu)\,\mbox{~~s.t.~~}\|f\|_C^0\leq 1,\, Lip(f)\leq 1}\label{e-dualgw}
\eqn
\ep
We recall that the $L^1$ distance satisfies a dual formula too, that is
\bqn
|\mu-\nu|=\sup\Pg{\int f d(\mu-\nu)\,\mbox{~~s.t.~~}\|f\|_C^0\leq 1}\label{e-dualL1}
\eqn

We also have this useful estimate to bound integrals. See \cite{gw}.
\bl Let $f\in C^0(\R^d)\cap\mathrm{Lip}(\R^d)$. It then holds
\bqn
\int f d(\mu-\nu)\leq \max\Pg{\|f\|_{C^0},\mathrm{Lip}(f)} \gw(\mu,\nu)\label{e-dual}.
\eqn
\el

We end this section by giving useful estimates both for the standard and generalized Wasserstein distances $W_p$ and $\gw$ under flow actions. Proofs are given in \cite{gw,gw2}.

\bp \label{p-flow} Let $v_t,w_t$ be two time-varying vector fields, uniformly Lipschitz with respect to the space variable, and $\phi^t,\psi^t$ the flow generated by $v,w$ respectively. Let $L$ be the Lipschitz constant of $v$ and $w$, i.e. $|v_t(x)-v_t(y)|\leq L |x-y|$ for all $t$, and similarly for $w$. Let $\mu,\nu\in\Mu$. We have the following estimates for the standard Wasserstein distance
\bi
\i $W_p\Pt{\phi^t\#\mu,\phi^t\#\nu}\leq e^{Lt}W_p\Pt{\mu,\nu}$,
\i $W_p\Pt{\mu,\phi^t\#\mu}\leq t\|v\|_{C^0} |\mu|$,
\i $W_p\Pt{\phi^t\#\mu,\psi^t\#\nu}\leq e^{Lt}W_p\Pt{\mu,\nu}+\frac{e^{Lt}-1}{L}\,|\mu|\,\sup_{\tau\in[0,t]}\|v_t-w_t\|_{C^0} $.
\ei
We have the following estimates for the generalized Wasserstein distance
\bi
\i $\gw(\phi^t\#\mu,\phi^t\#\nu)\leq e^{Lt}\gw(\mu,\nu)$,
\i $\gw(\mu,\phi^t\#\mu)\leq  t\|v\|_{C^0} |\mu|$,
\i $\gw(\phi^t\#\mu,\psi^t\#\nu)\leq e^{Lt}\gw(\mu,\nu)+\frac{e^{Lt}-1}{L}\,|\mu|\,\sup_{\tau\in[0,t]}\|v_t-w_t\|_{C^0}$.
\ei
\ep

\subsection{Measure Differential Equations with Probability Vector Fields} \label{s-MDE}

In this section, we summarize the main results and tools about PVFs, introduced in \cite{MDE}. We slightly enlarge the setting of \cite{MDE}, since we consider general measures with finite mass and not only probability measures.

We first recall the definition of a solution to the Cauchy problem
\bqn
\dot \mu=V[\mu],\qquad \mu(0)=\mu_0.
\eqnl{e-MDEcauchy}
\bdeff Fix a final time $T>0$. A solution to \r{e-MDEcauchy} is a map $\mu:[0,T]\to \M(\R^n)$ such that $\mu(0)=\mu_0$ and the following holds: for each $f\in C^\infty_c(\R^n)$
\bi
\i the integral $\int_{T\R^n} (\nabla f(x)\cdot v)\, dV[\mu(s)](x,v)$ is defined for almost every $s\in[0,T]$;
\i the map $s\to \int_{T\R^n} (\nabla f(x)\cdot v)\, dV[\mu(s)](x,v)$ belongs to $L^1([0,T])$;
\i the map $t \to \int f\, d\mu(t)$ is absolutely continuous, and it satisfies 
\bqn
\frac{d}{dt} \int_{\R^n} f\,d \mu(t)=\int_{T\R^n} (\nabla f(x)\cdot v)\, dV[\mu(t)](x,v)
\eqnl{e-weakPDE}
for almost every $t\in[0,T]$.
\ei
\edeff

We now recall the definition of the pseudo-distance $\W$, that will be useful in the following.
\bdeff \label{d-W} Let $V_1,V_2\in \M(T\R^n)$ with $|V_1|=|V_2|$. Denote by $\mu_1=\pi_1\#V_1$ and $\mu_2=\pi_1\#V_2$ the projection of the PVF on the base space. Define
\bqn
\W(V_1,V_2):=\inf\Pg{\int_{T\R^n\times T\R^n} |v-w|\, d p(x,v,y,w)\, \mbox{ such that }p\in P(V_1,V_2) \mbox{ and }\pi_{13}\#p\in P^{opt}(\mu_1,\mu_2)}.
\eqnn
\edeff
Clearly, such functional is not a distance, see examples in \cite{MDE}. Nevertheless, we will see in the following that the local Lipschitz condition {\bf (V2)} will ensure existence of solutions to \r{e-MDEcauchy}. Observe that it also holds
\bqn
W(V_1,V_2)\leq \W(V_1,V_2)+W(\pi_1\#V_1,\pi_1\#V_2).
\eqnl{e-WW}
See \cite{MDE} for more details.\\

We now address the problem of existence of solutions to \r{e-MDEcauchy}. The idea developed in \cite{MDE} is to define a semigroup of solutions as the limit of approximated ones. We first describe precisely the discretization method, that will be also useful in the following.
\bdeff\label{d-discrete} Fix $N\in\N$ and define the time step size $\Delta_N=\frac1N$, the velocity step size $\Delta^v_N=\frac1N$ and the space step size $\Delta^x_N=\Delta^v_N \Delta_N=\frac{1}{N^2}$. Define $x_i$ the $(2N^3+1)^n$ equispaced discretization points of $(\Z^n/(N^2))\cap [-N,N]^n$, and $v_j$ the $(2N^2+1)^n$ equispaced discretization points of $(\Z^n/N)\cap [-N,N]^n$.

Define $\M_N^x\subset \M(\R^n)$ the space of measures of $\R^n$ with support on the set of points $x_i$, and $\M_N^v\subset \M(\R^{2n})$ the space of measures of $\R^{2n}$ with support on the set of points $(x_i,v_j)$,

Define the discretization operator in the space variable $A^x_N:\M(\R^n)\to \M_N^x$ as follows
$$\A^x_N(\mu):=\sum_i m^x_i(\mu) \delta_{x_i},$$
where $m^x_i(\mu):=\mu(x_i+Q)$ with $Q=\left[0,\frac{1}{N^2}\right)^n$. Define the discretization operator in the velocity variable $A^v_N:\M(\R^{2n})\to \M_N^v$ as follows
\bqn
\A^v_N(V):=\sum_{i,j} m^v_{ij}(V) \delta_{(x_i,v_j)},
\eqnl{e-defAv}
where $m^v_{ij}(V):=V((x_i+Q)\times (v_j+Q'))$ with $Q'=\left[0,\frac{1}{N}\right)^n$.
\edeff
The first property of such discretization is that it introduces an arbitrarily small error in the Wasserstein distance.
\bp \label{p-errorediscr} Given $\mu\in \M_c(\R^n)$ and $V[\mu]\in \M_c(T\R^n)$, for a sufficiently large $N$ it holds
\bqn 
W(\mu,\A^x_N(\mu))\leq |\mu|\Delta^x_N,\qquad W(V[\mu],\A^v_N(V[\mu])\leq |\mu|\Delta^v_N.
\eqnn
\ep
\bproof The proof with $\mu$ and $V[\mu]$ being probability measures is given in \cite{MDE}. The generalization to measures with finite mass is straightforward.
\eproof

One can then define an approximated solution (called the Lattice Approximate Solution) to \eqref{e-MDEcauchy} via an explicit Euler scheme. 
\bdeff\label{d-LAS}
Given the Cauchy problem \r{e-MDEcauchy}, we define the following Lattice Approximate Solution $\mu^N$: we set $\mu^N(0):=\A^x_N(\mu_0)$, then recursively
$$\mu^N((k+1)\Delta_N):=\sum_{i,j} m^v_{ij}(V[\mu^N(k\Delta_N)])\delta_{x_i+\Delta_N v_j},$$
and for intermediate times $t\in[0,\Delta_N)$ we define
$$\mu^N(k\Delta_N+t):=\sum_{i,j} m^v_{ij}(V[\mu^N(k\Delta_N)])\delta_{x_i+t v_j}.$$
\edeff
We are now ready to state the existence of a solution to \r{e-MDEcauchy} as a limit of the Lattice Approximate Solutions introduced above.
\bt\label{t-exPVF} Let a PVF $V:\M(\R^n)\to \M(T\R^n)$ be given, satisfying {\bf (V)} where $\GW$ is replaced by $\W$. Then,  there exists a  Lipschitz semigroup of solutions to \r{e-MDEcauchy}, obtained as uniform-in-time limit of Lattice Approximate Solutions for the Wasserstein Metric.\et
\bproof The first key observation is that both $A^x_N$ and $A^v_N$ are operators preserving the mass for $N$ sufficiently large, i.e. $\A^x_N(\mu)(\R^n)=\mu(\R^n)$ and similarly for the PVF. As a consequence, the mass of $\mu^N(t)$ coincides with $\mu^N(0)$, that in turn coincides with $\mu_0$ for $N$ sufficiently large.

If $\mu_0(\R^n)=1$, then the whole sequence $\mu^N(t)$ is in $\P_c(\R^n)$, and one can apply the proof of \cite[Theorem 4.1]{MDE}. Otherwise, rescale the mass by defining $\nu^N(t)=\frac{1}{\mu_0(\R^N)}\mu^N(t)$, apply the previous case to define $\nu(t)$ and prove that $\mu(t)= {\mu_0(\R^N)} \nu(t)$ is a solution to \r{e-MDEcauchy}.\eproof

We now recall the definition of Dirac germs, that permits to address the problem of uniqueness of the solution to \r{e-MDEcauchy}. We also give the definition of semigroup compatible with the germ.
\bdeff \label{d-germ} Fix a PVF $V$. Define $\mathcal{M}^D:=\Pg{\mu\in \M(\R^n)\mbox{ such that } \mu=\sum_{l=1}^m m_l \delta_{x_l}}$ the space of measures composed of Dirac deltas. A Dirac germ $\gamma$ compatible with $V$ is a map assigning to each $\mu\in \mathcal{M}^D$ a Lipschitz curve $\gamma_\mu:[0,\eps(\mu)]\to \M(\R^n)$, with the following conditions:
\bi
\i $\eps(\mu)>0$ is uniformly positive for measures with uniformly bounded support;
\i $\gamma_\mu$ is a solution to \r{e-MDE}.
\ei
\edeff
\bdeff \label{d-compatible} Fix a PVF $V$ satisfying {\bf (V1)}, a final time $T>0$ and a Dirac germ $\gamma$. A semigroup for \r{e-MDE} is said to be compatible with $\gamma$ if one has the following property: for each $R,M>0$ there exists $C(R,M)$ such that the space $\M^D_{R,M}:=\Pg{\mu\in\M^D\mbox{~~s.t.~~}\supp(\mu)\in B(0,R), |\mu|\leq M}$ satisfies
\bqn
\mbox{for all }t\in[0,\inf_{\mu\in\M^D_{R,M}}\eps(\mu)] \mbox{ one has} \sup_{\mu\in\M^D_R} W(S_t\mu,\gamma_\mu(t))\leq C(R,M) t^2.
\eqnl{e-compat}
\edeff

We are now ready to prove the main result about uniqueness of solutions to \r{e-MDEcauchy}.
\bt Consider a PVF satisfying {(V1)} and fix a Dirac germ $\gamma$. There exists at most one Lipschitz semigroup $S_t$ of solutions to \r{e-MDEcauchy} compatible with $\gamma$.
\et
\bproof First observe that uniform boundedness of the support and the weak formulation \eqref{e-weakPDE} when choosing $f=1$ on $\cup_{t\in[0,T]} \supp(\mu(t))$ imply that the mass $\mu(t)(\R^n)$ is constant along trajectories of \r{e-MDE}. Thus, the Dirac germ satisfies conservation of mass too.

Apply now the proof of Theorem 5.1 in \cite{MDE} for an initial data being a probability measure, with the Dirac germ restricted to probability measures. For initial data with general finite mass, apply the rescaling trick described in the proof of Theorem \ref{t-exPVF} both to the initial data and the Dirac germ.
\eproof
\subsection{Measure Equations with sources}

In this section, we briefly study the measure equation with source 
\bqn
\dot \mu=s[\mu],\qquad \mu(0)=\mu_0.
\eqnl{e-sourcecauchy}

The goal is to prove that condition {\bf (s)} in Theorem \ref{t-existence} ensures existence and uniqueness of a solution to \r{e-sourcecauchy}. This is indeed a particular case of a more general result, stated in \cite{gw}, in which a transport term is added too. For our future use, we prove the statement with the same discretization method of Lattice Approximate Solution introduced in Definition \ref{d-discrete}.
\bp Fix $T>0$. Let the source $s:\M(\R^n)\to \M(\R^n)$ satisfy Hypotheses {\bf (s)} in Theorem \ref{t-existence}. Then, there exists a unique solution to \r{e-sourcecauchy}.

Moreover, such solution is the uniform-in-time Wasserstein limit for $N\to\infty$ of Lattice Approximate Solutions $\mu^N:[0,T]\to\M(\R^n)$ defined as follows: Define $\mu^N(0):=\A^x_N(\mu_0)=\sum_{i} m^x_i(\mu_0)\delta_{x_i}$, then recursively 
\bqn
\mu^N((k+1)\Delta_N)=\mu^N(k\Delta_N)+\Delta_N \A^x_N(s[\mu^N(k\Delta_N)]).
\eqnl{e-schemasource}
We also define the time-interpolated solution for $t\in[0,\Delta_N]$ as follows: $\mu^N(k\Delta_N+t)=\mu^N(k\Delta_N)+t \A^x_N(s[\mu(k\Delta_N)]).$
\ep
\bproof We first prove existence of a solution, based on the Lattice Approximate Solution. We prove that $\mu^N$ is a sequence of equi-Lipschitz and equi-bounded curves in $C^0([0,T],\M(X))$, where $X$ is a compact subset of $\R^n$ and the space $\M(X)$ is endowed with the generalized Wasserstein distance $\gw$. For $\tau,\sigma\in[0,\Delta_N]$ it holds
\bqn
&&\gw(\mu^N(k\Delta_N+\tau),\mu^N(k\Delta_N+\sigma))=\gw(\tau A^x_N(s[\mu^N(k\Delta_N)]),\sigma A^x_N(s[\mu^N(k\Delta_N)]))\leq\nn
&& |\tau-\sigma| \Pabs{A^x_N(s[\mu^N(k\Delta_N)])}\leq |\tau-\sigma| \Pabs{s[\mu^N(k\Delta_N)]}.
\eqnl{e-proofsource}
We are then left to prove that $\Pabs{s[\mu^N(k\Delta_N)]}$ is uniformly bounded for $k\Delta_N\in[0,T]$. It is sufficient to observe that \r{e-gwL1} and hypothesis {\bf (s1)}, together with \r{e-proofsource}, imply
\bqn
\Pabs{s[\mu^N((k+1)\Delta_N)]}\leq\Pabs{s[\mu^N(k\Delta_N)]}+\gw(s[\mu^N(k\Delta_N+\tau)],s[\mu^N(k\Delta_N+\sigma)]),\leq (1+L\Delta_N)\Pabs{s[\mu^N(k\Delta_N)]}.
\eqnn
hence recursively $\Pabs{s[\mu^N((k+1)\Delta_N)]}\leq e^{LT}\Pabs{s[\mu^N(0)]}\leq e^{LT}|s[\mu_0]|$.

We now prove that there exists $R'$ such that  $\supp(\mu^N(t))\subset B(0,R')$ for all $N$ and $t\in[0,T]$. Eventually enlarging the radius $R$ given in hypothesis {\bf (s2)}, one can assume that $\supp(\mu_0)\subset B(0,R)$. Thus, the approximation operator $\A^x_N$ satisfies $\supp(A^x_N(\mu_0))\subset B(0,R+1)$, as well as $\supp(A^x_N(s[\mu]))\subset B(0,R+1)$ for any $\mu\in \M(\R^n)$. Since sum of measures with the same support gives a measure with the same support, one can easily prove by induction that measures $\mu^N(k\Delta_N+\tau)$ defined by the scheme \r{e-schemasource} all have support contained in $B(0,R')$ with $R'=R+1$.

Choose now $X=\overline{B(0,R')}$, that is a compact space. Then, $\M(X)$ is complete when endowed with the generalized Wasserstein distance $\gw$, see Proposition \ref{p-complete}. The sequence $\mu^N$ is equi-Lipschitz in $\M(X)$, due to \r{e-proofsource}, and equi-bounded, since masses are equi-bounded. Then, there exists a converging subsequence, converging to some $\mu^*\in C^0([0,T],\M(X))$. Recall that such convergence with respect to $\gw$ coincides with weak convergence of measures.

We now prove that $\mu^*$ is a solution to \r{e-sourcecauchy}. We first observe that $\gw(\mu^N(0),\mu_0)\leq |\mu_0|\Delta^x_N$ for $N$ sufficiently large implies $\mu^*(0)=\mu_0$. We now prove that for all $f\in C^\infty_c(\R^n)$ and $\tau,\sigma\in[0,T]$ with $\sigma>\tau$, it holds
\bqn
\int f(x)\,d(\mu^*(\sigma)-\mu^*(\tau))-\int_\tau^\sigma dt\int f(x)\, ds[\mu^*(t)]=0
\eqnl{e-sourceweak}
The definition \r{e-schemasource} implies that, for $\sigma, \tau\in[0,\Delta_N]$, it holds
\bqn
\int f(x)\,d(\mu^N(k\Delta_N+\sigma)-\mu^*(k\Delta_N+\tau))=(\sigma-\tau) \int f(x)\, dA^x_N(s[\mu^N(k\Delta_N)]).
\eqnl{e-sourcewd}
We then have
\bqn
&&\Pabs{\int f(x)\,d(\mu^*(\sigma)-\mu^*(\tau))-\int_\tau^\sigma dt\int f(x)\, ds[\mu^*(t)]}\leq
\Pabs{\int f(x)\,d(\mu^*(\sigma)-\mu^N(\sigma))}+\nn
&&\Pabs{\int f(x)\,d(\mu^*(\tau)-\mu^N(\tau))}+\int_\tau^\sigma dt \Pabs{\int f(x)\,d(s[\mu^N(k^t_N\Delta_N)]-s[\mu^*(t)]) }+\nn
&&\Pabs{\int f(x)\,d(\mu^N(\sigma)-\mu^N(\tau))-\int_\tau^\sigma dt \int f(x)\,ds[\mu^N(k^t_N\Delta_N)]},
\eqnl{e-sourceweak2}
where $k^t_N$ is the largest integer such that $t\geq k\Delta_N$, i.e. $k^t_N=\lfloor \frac{t}{\Delta_N}\rfloor$. The first two terms converge to zero since $\mu^N\weak\mu^*$, while the last term is identically zero due to \r{e-sourcewd}. For the third term, observe that it holds
\bqn
\Pabs{\int f(x)\,d(s[\mu^N(k^t_N\Delta_N)]-s[\mu^*(t)]) }\leq \|f\|_{C^1} L\Pt{\gw(\mu^N(k^t_N\Delta_N),\mu^N(t))+\gw(\mu^N(t),\mu^*(t)) },\eqnl{e-sss}
where we used condition {\bf (s1)} about the Lipschitz continuity of $s$, as well as the dual formulation \r{e-dual} for $\gw$. The proof now follows from observing that both the terms in the right hand side of \r{e-sss} converge to zero: the first satisfies 
\bqn
\gw(\mu^N(k^t_N\Delta_N),\mu^N(t))\leq K' |t-k^t_N\Delta_N|\leq K'\Delta_N\to 0
\eqnn
for the constant $K'=e^{LT}|s[\mu_0]|$ given by \r{e-proofsource}. The second converges to zero since $\gw$ metrizes weak convergence.

We now prove uniqueness of the solution, by proving continuous dependence on the initial data for \r{e-sourcecauchy}. Consider two solutions $\mu(t),\nu(t)$ to  \r{e-sourcecauchy} with initial data $\mu_0,\nu_0$ respectively. Using the the weak formulation of \r{e-sourcecauchy}, it holds
\bqn
\int f(x)\,d(\mu(t)-\nu(t))=\int f(x)\,d(\mu_0-\nu_0)+\int_0^t d\tau f(x) \,d(s[\mu(\tau)]-s[\nu(\tau)]).
\eqnl{e-ddual}
Choose now a sequence $f_n$ with $\|f_n\|_C^0\leq 1,\, Lip(f_n)\leq 1$ realizing $\gw(\mu(t),\nu(t))$ in its dual formulation \r{e-dualgw}. Then, equation \r{e-ddual} reads as $\gw(\mu(t),\nu(t))\leq \gw(\mu_0,\nu_0)+\int_0^t L \gw(\mu(\tau),\nu(\tau))$, where we used Lipschitz continuity of $s$. Since both $\mu(t),\nu(t)$ are Lipschitz with respect to time, a direct application of the Gronwall lemma implies $\gw(\mu(t),\nu(t))\leq e^{Lt}\gw(\mu_0,\nu_0)$, that in turn implies uniqueness of the solution to  \r{e-sourcecauchy}.
\eproof

\subsection{The operator $\GW$}\label{s-GW}

In this section, we define the operator $\GW$. The Lipschitz condition {\bf (V2)} with respect to such operator will be crucial to ensure existence of a solution to \r{e-MDES}.
Then we can define:
\bdeff Let $V_1,V_2\in \M(T\R^n)$. Denote by $\mu_1=\pi_1\#V_1$ and $\mu_2=\pi_1\#V_2$ the projection of the PVF on the base space. For each pair $(\tilde V_1,\tilde V_2)$ satisfying $\tilde V_1\leq V_1$ and $\tilde V_2\leq V_2$, denote by $\tilde\mu_i=\pi_1\#\tilde V_i$. Define
\bqn
\GW(V_1,V_2)&:=&\inf\left\{\int_{T\R^n\times T\R^n} |v-w|\, d p(x,v,y,w)\, \mbox{ such that }\tilde V_1\leq V_1,\ \tilde V_2\leq V_2,
p\in P(\tilde V_1,\tilde V_2),
\right.\label{e-GW}\\
&& \left. \gw(\mu_1,\mu_2)=|\mu_1-\tilde \mu_1|+W(\tilde\mu_1,\tilde\mu_2)+|\mu_2-\tilde\mu_2| \mbox{ and }\pi_{13}\#p\in P^{opt}(\tilde\mu_1,\tilde\mu_2)\right\}.
\eqnn
\edeff

In the definition above, one can observe that the condition $\gw(\mu_1,\mu_2)=|\mu_1-\tilde \mu_1|+W(\tilde\mu_1,\tilde\mu_2)+|\mu_2-\tilde\mu_2|$ is equivalent to state that $\tilde\mu_1,\tilde\mu_2$ is a minimizer in Definition \ref{d-gw}. 

\brem One might require the minimization of the functional $|V_1-\tilde V_1|+\int_{T\R^n\times T\R^n} |v-w|\, d p(x,v,y,w)+|V_2-\tilde V_2|$ in \r{e-GW}, that seems more close to the definition of $\gw$. Nevertheless, recall that that $|V_i-\tilde V_i|=|\mu_i-\tilde \mu_i|$. As a consequence, when the choice of the minimizer for $\gw(\mu_1,\mu_2)$ is unique, there is no difference for the minimization of the two functionals. When minimizers for $\gw(\mu_1,\mu_2)$ are not unique, this would introduce two additional terms $|V_i-\tilde V_i|$ in the right hand side of \r{e-WW2}, thus providing a less restricitve inequality.

Moreover, the chosen definition of $\GW$ is correct to prove existence of a solution to \r{e-MDES}, that is the main goal of this paper. This definition would indeed be natural in the estimate \r{e-Prop2} below.
\erem

When two measures $\mu_1,\mu_2$ have the same mass $|\mu_1|=|\mu_2|$ and have sufficiently near supports, one can choose $\tilde\mu_i=\mu_i$ among the minimizers of $\gw(\mu_1,\mu_2)$. Moreover, if $\mu_1\perp\mu_2$, i.e. the two measures have no shared mass, such choice is the unique minimizer. In this case, the choice $\tilde V_i=V_i$ is unique too, and the operator $\GW$ coincides with $\W$.

Similarly to $\W$, the operator $\GW$ is then not a distance: the same counterexamples with sufficiently near supports and no shared mass can be found. For example, choose $\mu_1=\delta_0, \mu_2=\delta_\epsilon, V_1=\delta_0\otimes\delta_0, V_2=\delta_\epsilon\otimes \delta_0$, with $\epsilon>0$ sufficiently small. The unique minimizer in \r{e-GW} is given by $\tilde\mu_i=\mu_i$ and $\tilde V_i=V_i$, that in turn gives $\GW(V_1,V_2)=0$, even though $V_1\neq  V_2$.

Again, similarly to the estimate \r{e-WW} between the standard generalized Wasserstein and the operator $\W$, one can easily prove the following proposition.
\bp Given $V_1,V_2\in\M(T\R^n)$ two PVFs, it holds
\bqn
\gw(V_1,V_2)\leq \GW(V_1,V_2)+\gw(\pi_1\#V_1,\pi_1\#V_2).
\eqnl{e-WW2}
\ep
\bproof First estimate
$\gw(V_1,V_2)$ from above by choosing $\tilde V_1,\tilde V_2$ minimizers of $\GW(V_1,V_2)$ in the right-hand side of \r{e-gw}. Similarly, estimate the Wasserstein distance $W(\tilde V_1,\tilde V_2)$ from above by using the transference plan $p$ realizing $\GW(V_1,V_2)$. It then holds
\bqn
\gw(V_1,V_2)&\leq& |V_1-\tilde V_1|+W(\tilde V_1,\tilde V_2)+|V_2-\tilde V_2|\leq |\mu_1-\tilde \mu_1|+ \W(\tilde V_1,\tilde V_2)+W(\tilde\mu_1,\tilde\mu_2)+ |\mu_2-\tilde \mu_2|=\nn
&=&\gw(\mu_1,\mu_2)+\W(\tilde V_1,\tilde V_2).
\eqnn
We used here that $|V_i-\tilde V_i|=|\mu_i-\tilde \mu_i|$ since $\tilde V_i\leq V_i$ and $\tilde\mu_i\leq \mu_i$. We also used \r{e-WW} and the fact that $\tilde\mu_1,\tilde\mu_2$ is a minimizer for $\gw(\mu_1,\mu_2)$.

We are now left to prove that $\W(\tilde V_1,\tilde V_2)=\GW(\tilde V_1,\tilde V_2)$. It clearly holds 
$\W(\tilde V_1,\tilde V_2)\geq\GW(\tilde V_1,\tilde V_2)$, since the minimization on the right hand side takes place in a larger space. By contradiction, if a strict inequality holds true, there exists a decomposition $(\hat V_1,\hat V_2)$ satisfying $\hat V_i < \tilde V_i$ and minimizing $\GW(\tilde V_1,\tilde V_2)$. This means that there exists $q\in P(\hat V_1,\hat V_2)$ such that $$\int_{T\R^n\times T\R^n} |v-w|\, d q(x,v,y,w)<\W(\tilde V_1,\tilde V_2)$$ and, by defining $\hat\mu_i=\pi_1\# \hat V_i$, it holds
$$\gw(\tilde\mu_1,\tilde\mu_2)=|\tilde\mu_1-\hat\mu_1|+W(\hat\mu_1,\hat\mu_2)+|\tilde\mu_2-\hat\mu_2|,$$ with $\pi_{13}\#q\in P^{opt}(\hat\mu_1,\hat\mu_2)$. Also recall that $\gw(\tilde\mu_1,\tilde\mu_2)\leq W(\tilde\mu_1,\tilde\mu_2)$, by \r{e-gww}. 

Observe now that $\hat V_1,\hat V_2$ is a possible decomposition to estimate $\GW(V_1,V_2)$ in \r{e-GW}, with transference plan $q$. Indeed, it first holds $\hat V_i\leq \tilde V_i\leq V_i$, thus $|\mu_i-\hat\mu_i|=|\mu_i-\tilde\mu_i|+|\tilde\mu_i-\hat\mu_i|$. This implies that
\bqn
&&|\mu_1-\hat \mu_1|+ W(\hat\mu_1,\hat\mu_2)+ |\mu_2-\hat \mu_2|=\nn
&&|\mu_1-\tilde\mu_1|+|\tilde\mu_1-\hat\mu_1|+W(\hat\mu_1,\hat\mu_2)+|\tilde\mu_2-\hat\mu_2|+|\mu_2-\tilde\mu_2|=\nn
&&|\mu_1-\tilde\mu_1|+\gw(\tilde\mu_1,\tilde\mu_2)+|\mu_2-\tilde\mu_2|\leq 
|\mu_1-\tilde\mu_1|+W(\tilde\mu_1,\tilde\mu_2)+|\mu_2-\tilde\mu_2|=\gw(\mu_1,\mu_2),
\eqnn
i.e.  the decomposition $\hat\mu_1,\hat\mu_2$ realizes the minimizer of $\gw(\mu_1,\mu_2)$. Since $\pi_{13}\#q\in P^{opt}(\hat\mu_1,\hat\mu_2)$, one can write by the contradiction hypothesis that it holds 
$$\int_{T\R^n\times T\R^n} |v-w|\, d q(x,v,y,w)=\GW(\tilde V_1,\tilde V_2)<\W(\tilde V_1,\tilde V_2)=\GW(V_1,V_2).$$
This contradicts the definition of $\GW$ as the infimum of the functional $\int_{T\R^n\times T\R^n} |v-w|\, d p(x,v,y,w)$.
\eproof

\section{Proof of the main theorems} \label{s-proof}
In this section, we prove the main results of this article, that are Theorem \ref{t-existence} about existence of solutions to \r{e-MDES} and Theorem \ref{t-uniqueness} about uniqueness.

\subsection{Existence -  Proof of Theorem \ref{t-existence}} \label{s-existence}
In this section, we prove Theorem \ref{t-existence} about existence of solutions to \r{e-MDES}. The idea is to define Lattice Approximate Solutions described in Definition \ref{d-LAS}, then pass to the limit. This procedure proved useful for each term in \r{e-MDES} separately, namely the PVF studied in \r{e-MDE} and the source in \r{e-source}.

{\it Proof of Theorem \ref{t-existence}}.  We first fix an initial data $\mu_0$ and prove the existence of a solution to \r{e-MDES} with initial data $\mu_0$, that has bounded support and is Lipschitz with respect to time. This corresponds to prove Properties 1-2-3a-3c in the Definition \ref{d-semigroup} of semigroups for \r{e-MDES}. We will then prove Property 3b.

Fix an initial data $\mu_0$. For each $N$, define the following approximated solution $\mu^N$, based on the discretization in Definition \ref{d-discrete}: set $\mu^N(0):=\A^x_N(\mu_0)$, then recursively
\bqn
\mu^N((k+1)\Delta_N):=\sum_{i,j} m^v_{ij}(V[\mu^N(k\Delta_N)])\delta_{x_i+\Delta_N v_j}+\Delta_N \A^x_N(s[\mu^N(k\Delta_N)]).
\eqnl{e-schema}
Also define the interpolated measure for $\tau\in [0,\Delta_N]$ as follows:
\bqn\mu^N(k\Delta_N+\tau):=\sum_{i,j} m^v_{ij}(V[\mu^N(k\Delta_N)])\delta_{x_i+\tau v_j}+\tau \A^x_N(s[\mu^N(k\Delta_N)]).
\eqnl{e-schemat}
Clearly, the first term on the right hand side corresponds to the transport by the PVF $V$, while the second term corresponds to the source term $s$. We now prove that the sequence $\mu^N(t)$ is equi-bounded and equi-integrable on a complete space, to apply the Ascoli-Arzel\`a theorem.

We first prove that, for a fixed $T>0$, the measures $\mu^N(t)$ are all supported in a compact set. Choose the radius $R$ in hypothesis {\bf (s2)} giving the maximal support of $s[\mu]$. Then, eventually enlarging $R$, one can assume that $\supp(\mu(0))\subset B(0,R)$. This implies $\supp(\mu^N(0))\subset B(0,R+\sqrt{\frac{n}2 \Delta^x_N})\subset B(0,R+1)$, where $n$ is the dimension of the space $\R^n$ and $N$ is chosen sufficiently large.  Observe the following simple estimate: if $\supp(\mu^N(k\Delta_N))\subset B(0,r)$ with $r>R+1$, then it holds $\supp(\mu^N(k\Delta_N+\tau))\subset B(0,r+\Delta_N C(1+r))$ for $\tau\in(0,\Delta_N)$. Indeed:
\bi
\i for each term $i,j$ in  the first term it holds that $(x_i,v_j)\in\supp(V[\mu^N(k\Delta_N)])$ implies $|v_j|\leq C(1+r)$ by Hypothesis {\bf (V1)}, hence $\supp(\delta_{x_i+\tau v_j})\subset B(0,r+\Delta_N C(1+r))$;
\i for the second term it holds $\supp(s[\mu^N(k\Delta_N)])\subset B(0,R)$, hence $\Delta_N \A^x_N(s[\mu^N(k\Delta_N)])\subset B(0,R+1)\subset B(0,r)$.
\ei
Since summing measures with the same support is a closed operation, it  holds $\supp(\mu^N(k\Delta_N+\tau))\subset B(0,r+\Delta_N C(1+r))$. Eventually replacing $r$ by $\max\Pg{1,r}$, it holds by induction $\supp(\mu^N(k\Delta_N+\tau))\subset B(0, r(1+2C\Delta_N)^{k+1})$. Since $k\leq \frac{T}{\Delta_N}+1$, this implies $\supp(\mu^N(t))\subset B(0, r (1+2C)^2 e^{2CT})$ for all $t\in [0,T]$. Observe that the space $X=\overline{B(0,r(1+2C)^2e^{2CT})}$ is compact. Then, the space $M(X)$ of measures with finite mass endowed with the generalized Wasserstein distance $\gw$ is complete, see Proposition \ref{p-complete}. Moreover, if we prove that the limit of a subsequence of $\mu^N$ exists, then it satisfies {\bf Property 3a} in Definition \ref{d-semigroup}.

We now prove that the sequences $\mu^N(t)$ are equi-Lipschitz in time with respect to the distance $\gw$. We also prove that the masses $|\mu^N(t)|$ are uniformly bounded. First observe that the operator $\mu\to \A^x_N(\mu)$ does not increase the mass of $\mu$. The same property holds for the operator $\mu\to\sum_{i,j}m_{ij}^v (V[\mu])\delta_{x_i+\tau v_j}$.  Then, by the explicit expressions \r{e-schema}-\r{e-schemat} for $\mu^N(t)$, it holds 
\bqn
&&\gw(\mu^N(k\Delta_N+\tau),\mu^N(k\Delta_N+\sigma))\leq \sum_{i,j} m_{ij}^v(V[\mu^N(k\Delta_N)]) |\tau-\sigma| |v_j| + |\tau-\sigma| |s[\mu^N(k\Delta_N)]|\leq \nn
&&\leq|\tau-\sigma|\,|\mu^N(k\Delta_N)|\,\sup_{j} |v_j|+ |\tau-\sigma| ( \gw(s[\mu^N(k\Delta_N)],s[\mu_0])+|s[\mu_0]|),
\eqnn
for $\tau,\sigma\in[0,\Delta_N]$. Here, we estimated the generalized Wasserstein distance by decomposing it in the Wasserstein distance for the transport term given by the PVF $V$, and the $L^1$ distance for the source term given by $s$. Use now \r{e-gwL1} to estimate $|\mu^N(k\Delta_N)|\leq |\mu_0|+\gw(\mu^N(k\Delta_N),\mu_0)$. Use uniform boundedness of the supports for $\mu^N(k\Delta_N)$ and Hypothesis {\bf (V1)} to estimate $\sup_{j} |v_j|\leq C_1:=C(1+\diam(X))$. Also use Hypothesis {\bf (s1)} to estimate $\gw(s[\mu^N(k\Delta_N)],s[\mu_0])\leq L\gw(\mu^N(k\Delta_N),\mu_0)$. This gives
\bqn
&&\gw(\mu^N(k\Delta_N+\tau),\mu^N(k\Delta_N+\sigma))\leq |\tau-\sigma| C_1(|\mu_0|+\gw(\mu^N(k\Delta_N),\mu_0))+\nn
&&|\tau-\sigma| (|s[\mu_0]|+ L\gw(\mu^N(k\Delta_N),\mu_0))=|\tau-\sigma|(C_2\gw(\mu^N(k\Delta_N),\mu_0)+C_2),
\eqnl{e-gwkp1}
where we choose $C_2$ a constant depending on $|\mu_0|$ and $X$ only, i.e. independent on $N$.

We now prove that $\gw(\mu^N(k\Delta_N),\mu_0)$ is bounded, uniformly in $N,k$. Observe that the sequence $\mu^N(0)$ for $N$ sufficiently large satisfies $\gw(\mu^N(0),\mu_0)\leq W(\mu^N(0),\mu_0)\leq |\mu_0|\Delta_N$, as a consequence of \r{e-gww} and Proposition \ref{p-errorediscr}. Thus, there exists a constant $C_3$ such that  $\gw(\mu^N(0),\mu_0)\leq C_3$ for all $N$. We now prove the estimate
\bqn
\gw(\mu^N(k\Delta_N),\mu_0)\leq (1+C_2\Delta_N)^k C_3+((1+C_2\Delta_N)^k-1),
\eqnl{e-munmu0}
by induction on $k$. The case $k=0$ is already proved. We now prove that, if the estimate holds for $k$, then it holds for $k+1$ too. Use \r{e-gwkp1} with $\tau=0, \sigma=\Delta_N$, that gives
\bqn
&&\gw(\mu^N((k+1)\Delta_N),\mu_0)\leq \gw(\mu^N((k+1)\Delta_N),\mu^N(k\Delta_N))+\gw(\mu^N(k\Delta_N),\mu_0)\leq\nn
&& \leq (1+C_2\Delta_N)\gw(\mu^N(k\Delta_N),\mu_0)+C_2\Delta_N\leq (1+C_2\Delta_N)^{k+1} C_3+(1+C_2\Delta_N) ((1+C_2\Delta_N)^k-1)+\nn
&&C_2\Delta_N= (1+C_2\Delta_N)^{k+1} C_3+((1+C_2\Delta_N)^{k+1}-1).
\eqnn
The estimate \r{e-munmu0} is now proved. Since $k\leq T/\Delta_N$, it also holds $\gw(\mu^N(k\Delta_N),\mu_0)\leq C_4:=e^{C_2T}C_3+(e^{C_2T}-1)$. Then, again by \r{e-gwkp1} and triangular inequalities, it holds
\bqn
\gw(\mu^N(t),\mu^N(s))\leq |t-s|C_5,
\eqnl{e-ul}
with $C_5:=C_2C_4+C_2$, i.e. uniform Lipschitz continuity with respect to $t$ of the sequence $\mu^N$. Moreover, \r{e-gwL1} also implies $|\mu^N(t)|\leq C_6:=|\mu_0|+TC_5$, hence uniform boundedness of the mass. As a consequence, Ascoli-Arzel\`a theorem implies existence of a converging subsequence $\mu^N$, that we do not relabel. Such limit $\mu^*(t)$ clearly satisfies {\bf Property 1} and {\bf Property 3c} in Definition \ref{d-semigroup}.

We now prove {\bf Property 2}, i.e. the fact that the limit is a solution to \r{e-MDES}. Since the limit is uniformly Lipschitz and with bounded mass, it is easy to prove that the two first properties in Definition \ref{d-sol} are satisfied, as well as the fact that the function $t\to\int_{\R^n} fd\mu(t)$ is absolutely continuous (and even Lipschitz) for each $f\in C^\infty_c(\R^n)$. We are left to prove that the limit $\mu^*$ satisfies \r{e-weak} for each $f\in C^\infty_c(\R^n)$ and almost each $t\in[0,T]$. For a fixed $f\in C^\infty_c(\R^n)$, define the operator $F^N$ for $\tau,\sigma\in[0,T]$ as follows:
\bqn
F^N(\tau,\sigma):=\int f (x)\,d(\mu^N(\sigma)-\mu^N(\tau)) -\int_\tau^\sigma dt\Pt{\int \nabla f(x)\cdot v \,dV[\mu^N(t)])+\int f(x)\,ds[\mu^N(t)]}.
\eqnl{e-pezzo}

For each $\tau,\sigma\in[k\Delta_N,(k+1)\Delta_N]$, the first term in the right hand side of \r{e-pezzo} coincides with
\bqn
&&\int f(x) \Pt{\sum_{i,j} m_{ij}^v(V[\mu^N(k\Delta_N)]) d(\delta_{x_i+(\sigma-k\Delta_N) v_j}-\delta_{x_i+(\tau-k\Delta_N) v_j})+ (\sigma-\tau) \,d\A^x_N(s[\mu^N(k\Delta_N)])}.
\eqnl{e-eccola}

Define
\bqn
g_{ij}(\alpha):=f(x_i+\alpha (\sigma-k\Delta_N) v_j)-f(x_i+\alpha (\tau-k\Delta_N) v_j),
\eqnl{e-g} that satisfies $g_{ij}(1)=\int f(x) d(\delta_{x_i+(\sigma-k\Delta_N) v_j}-\delta_{x_i+(\tau-k\Delta_N) v_j})$. By applying Taylor's theorem with Lagrange remainder to each $g_{ij}$, there exist $\alpha_{ij}\in(0,1)$ such that \r{e-g} coincides with
\bqn 
g_{ij}(1)&=&(\sigma-\tau) \nabla f(x_i)\cdot v_j + v_j\cdot ((\sigma-k\Delta_N)^2 Hf(x_i+\alpha_{ij}\sigma v_j)-(\tau-k\Delta_N)^2 Hf(x_i+\alpha_{ij}\tau v_j))\cdot v_j,
\eqnn
where $Hf$ is the Hessian of $f$.

We now estimate $F^N(\tau,\sigma)$ with $\tau,\sigma\in[k\Delta_N,(k+1)\Delta_N]$. It holds
\bqn
|F^N(\tau,\sigma)|&\leq& \Pabs{\sum_{i,j} g_{ij}(1) m_{ij}^v(V[\mu^N(k\Delta_N)])-\int_\tau^\sigma dt \int \nabla f(x)\cdot v \,dV[\mu^N(t)])}+\nn
&&\Pabs{ (\sigma-\tau)\int f(x) \,d\A^x_N(s[\mu^N(k\Delta_N)])-\int_\tau^\sigma dt \int f(x)\,ds[\mu^N(t)]}.
\eqnl{e-g2}
The first term on the right hand side of \r{e-g2} is bounded from above by $\int_\tau^\sigma I_1(t)\,dt+I_2(\tau,\sigma)$, where
\bqn
I_1(t)&:=&\Pabs{\sum_{i,j} \nabla f(x_i)\cdot v_j\,m_{ij}^v(V[\mu^N(k\Delta_N)])- \int \nabla f(x)\cdot v \,dV[\mu^N(t)])}=\nn
&=&\Pabs{\int \nabla f(x)\cdot  v\, d\A^v_N(V[\mu^N(k\Delta_N)])-\int \nabla f(x)\cdot v \,dV[\mu^N(t)])},\nn
I_2(\tau,\sigma)&:=&\Pabs{\sum_{i,j}  v_j\cdot ((\sigma-k\Delta_N)^2 Hf(x_i+\alpha_{ij}\sigma v_j)-(\tau-k\Delta_N)^2 Hf(x_i+\alpha_{ij}\tau v_j))\cdot v_j\, m_{ij}^v(V[\mu^N(k\Delta_N)])}
\eqnn
The second term on the right hand side of \r{e-g2} is bounded from above by $\int_\tau^\sigma I_3(t)\,dt$, where
\bqn
&&I_3(t):=\Pabs{\int f(x)  \,d\A^x_N(s[\mu^N(k\Delta_N)])- \int f(x)\,ds[\mu^N(t)]}.
\eqnn

By the duality formulas for the generalized Wasserstein and the $L^1$ distance \r{e-dual}-\r{e-dualL1}, the Lipschitz condition \r{e-ul} and the estimate $|t-k\Delta_N|\leq \Delta_N$, for a sufficiently large $N$ it holds 
\bqn
&& I_1(t)\leq \|\nabla f\cdot v\|_{C^1} \gw(\A^v_N(V[\mu^N(k\Delta_N)]),V[\mu^N(t)])\leq\nn
&&( \|f\|_{C^2}\|v\|_{C^1} (\gw(\A^v_N(V[\mu^N(k\Delta_N)]),V[\mu^N(k\Delta_N)])+\gw(V[\mu^N(k\Delta_N)],V[\mu^N(t)]))\leq\nn
&&\|f\|_{C^2} (C_1+1)(\Delta_N |V[\mu^N(k\Delta_N)]|+|t-k\Delta_N|  C_5)\leq \Delta_N\|f\|_{C^2} (C_1+1)(C_6+C_5).
\eqnl{e-I1}
We recall that the support of the velocities is bounded ($\sup_{j} |v_j|\leq C_1$), and that the projection condition $\pi\#V[\mu]=\mu$ implies $|V[\mu]|=|\mu|$. Similarly, for $I_2(\tau,\sigma)$, by removing and adding the term $(\tau-k\Delta_N)^2 Hf(x_i+\alpha_{ij}\sigma v_j))$, it holds 
\bqn
I_2(\tau,\sigma)&\leq& \sum_{i,j} \|f\|_{C^2} (nC_1)^2 ((\sigma-k\Delta_N)^2-(\tau-k\Delta_N)^2)m_{ij}^v(V[\mu^N(k\Delta_N)])+\nn
&&\sum_{i,j}(nC_1)^2(\tau-k\Delta_N)^2 \| Hf(x_i+\alpha_{ij}\sigma v_j)-Hf(x_i+\alpha_{ij}\tau v_j)\|_{\R^n,\R^n} m_{ij}^v(V[\mu^N(k\Delta_N)]),
\eqnn
where $\|\cdot\|_{\R^n,\R^n}$ is the operator norm. Apply the mean-value theorem to $Hf$, to have
\bqn
\| Hf(x_i+\alpha_{ij}\sigma v_j)-Hf(x_i+\alpha_{ij}\tau v_j)\|_{\R^n,\R^n}\leq\alpha_{ij}(\sigma-\tau)v_j\|f\|_{C^3}. 
\eqnn Recall that $\sigma-k\Delta_N,\tau-k\Delta_N\leq \Delta_N$ as well as $\alpha_{ij}\in(0,1)$. It then holds
\bqn
I_2(\tau,\sigma)&\leq& \|f\|_{C^2}(nC_1)^2 (\sigma-\tau)2\Delta_N |V[\mu^N(k\Delta_N)]|+(nC_1)^3 \Delta_N^2(\sigma-\tau)\|f\|_{C^3}  |V[\mu^N(k\Delta_N)]|\nn
&\leq& (\sigma-\tau)\Delta_N (nC_1)^3C_6 \|f\|_{C^3} (2+\Delta_N).\label{e-I2}
\eqn

Finally, for $I_3(t)$ it holds
\bqn
I_3(t)&\leq& \|f\|_{C^1}  \gw(\A^x_N(s[\mu^N(k\Delta_N)]),s[\mu^N(t)])\leq \nn
 && \|f\|_{C^1}(W(\A^x_N(s[\mu^N(k\Delta_N)]),s[\mu^N(k\Delta_N)])+\gw(s[\mu^N(k\Delta_N)],s[\mu^N(t)]))\leq\nn
 && \|f\|_{C^1}(\Delta_N|s[\mu^N(k\Delta_N)]|+LC_5\Delta_N).
\eqnl{e-I3}
Observe that \r{e-gwL1} implies
\bqn
\Pabs{s[\mu^N(k\Delta_N)]}&\leq& |s[\mu_0]|+\gw(s[\mu^N(k\Delta_N)]|,s[\mu_0])\leq |s[\mu_0]|+\nn
&& L(\gw(\mu^N(k\Delta_N),\mu^N(0))+\gw(\mu^N(0),\mu_0))\leq |s[\mu_0]|+LC_5T+L|\mu_0|\Delta_N.
\eqnl{e-I3b}

Merging \r{e-I1}-\r{e-I2}-\r{e-I3}-\r{e-I3b}, it then holds $$|F^N(\tau,\sigma)|\leq|\tau-\sigma|\Delta_N\|f\|_{C^3} C_7$$ for a suitable constant $C_7$. Take now a general pair $\tau,\sigma\in[0,T]$. For simplicity, assume that $\tau<\sigma$ and that $N$ is sufficiently large to have 
$$(k_1-1)\Delta_N\leq \tau< k_1\Delta_N\leq k_2\Delta_N< \sigma\leq(k_2+1)\Delta_N$$ for some $k_1,k_2\in\N\setminus\Pg{0}$. Our aim is to prove that it holds $\lim_{\sigma\to\tau} \mathcal{L}=0$, with 
\bqn
&&\mathcal{L}:=\frac{1}{\sigma-\tau}\Pt{\int f(x)\,d(\mu^*(\sigma)-\mu^*(\tau))-\int_\tau^\sigma dt\int\nabla f(x)\cdot v \,dV[\mu^*(t)]-\int_\tau^\sigma dt\int f(x)\,ds[\mu^*(t)]}.
\eqnn
Observe that, restricting ourselves to the converging subsequence $\mu^N\weak\mu^*$,  it holds 
\bqn
&&\mathcal{L}=\lim_{N\to\infty}\frac{1}{\sigma-\tau}\Pt{\int f(x)\,d(\mu^N(\sigma)-\mu^N(\tau))-\int_\tau^\sigma dt\int\nabla f(x)\cdot v \,dV[\mu^N(t)]-\int_\tau^\sigma dt\int f(x)\,ds[\mu^N(t)]}.
\eqnn
We used here continuity of both $V$ and $s$ with respect to weak convergence of measures, as a consequence of \r{e-WW2} and {\bf (V2)} for $V$, and {\bf (s1)} for $s$. Observe that $\mathcal{L}$ coincides with
\bqn
\lim_{N\to\infty}\frac{1}{\sigma-\tau}\Pt{F^N(\tau,k_1\Delta_N)+\sum_{k=k_1}^{k_2-1}F^N(k\Delta_N,(k+1)\Delta_N)+F^N(k_2\Delta_N,\sigma)},
\eqnn
thus
\bqn
\lim_{\sigma\to\tau} |\mathcal{L}|\leq 
\lim_{\sigma\to\tau} \frac{1}{|\sigma-\tau|}\lim_{N\to\infty} (|\tau-k_1\Delta_N|+\ldots+ |k_2\Delta_N-\sigma|)\Delta_N\|f\|_{C^3} C_7=0.
\eqnn

We finally prove {\bf Property 3b} in Definition \ref{d-semigroup}. Take two different data $\mu_0,\nu_0$ and build the Lattice Approximate Solutions $\mu^N,\nu^N$ according with scheme \r{e-schema}-\r{e-schemat}. Assume to have $N$ sufficiently large so that it holds 
\bqn
\A^x_N[\mu^N(k\Delta_N)]=\mu^N(k\Delta_N)
\eqnl{e-Ax}
for all $k\in \N$ with $k\leq T/\Delta_N$, and similarly for $\nu^N$. Such $N$ exists, for two reasons: first, $\mu^N,\nu^N$ have uniformly bounded support (Property 3a proved above), thus for $N$ sufficiently large $\A^x_N$ conserves the mass. Second, observe that $\mu^N(0)$ satisfies \r{e-Ax} by construction, and that if $\mu^N(k\Delta_N)$ satisfies it, then $\mu^N((k+1)\Delta_N)$ satisfies it too. Then, by induction, this holds for all $k$.

Similarly, we assume to have $N$ sufficiently large to have
\bqn
\pi_1\#\A^v_N(V[\mu^N(k\Delta_N)])=\pi_1\#V[\mu^N(k\Delta_N)],
\eqnl{e-Av}
and similarly for $\nu^N$. This is first based on the fact that uniform bounded supports of the measure $\mu^N_t$ (Property 3a proved above), together with support sublinearity {\bf (V1)} of $V$ implies uniform boundedness of the supports of the PVF $V[\mu^N_t]$, thus $\A^v_N$ conserves the mass for $N$ sufficiently large. As soon as $|\A^v_N(V[\mu^N(k\Delta_N)])|=|\pi_1\#V[\mu^N(k\Delta_N)]|$, one has that the support of $\pi_1\#\A^v_N(V[\mu^N(k\Delta_N)])$ coincides with $\A^x_N(\mu^N(k\Delta_N))$, since $\pi_1\#V[\mu]=\mu$ and the discretization \eqref{e-defAv} has the same effect of $\A^x_N$ on the base space. Then, \r{e-Ax} implies \r{e-Av}.

Moreover, since \r{e-Av} holds, one has that any transference plan $p \in P (\A^v_N(V[\mu^N(k\Delta_N)]),V[\mu^N(k\Delta_N)])$ that satisfies $\pi_{13}\#p\in P^{opt}(\mu^N(k\Delta_N),\mu^N(k\Delta_N))$ is indeed a.e.-$\mu^N(k\Delta_N)$ concentrated on the diagonal $\{x=y\}$ of $\R^n\times\R^n$. As a consequence, one can estimate $\GW(\A^v_N(V[\mu^N(k\Delta_N)]),V[\mu^N(k\Delta_N)])$ by choosing in the right hand side of \r{e-GW} the decomposition $\tilde V_1=\A^v_N(V[\mu^N(k\Delta_N)])$, $\tilde V_2=V[\mu^N(k\Delta_N)]$ and the transference plan $p$ realizing $W(V[\mu^N(k\Delta_N)],\A^v_N(V[\mu^N(k\Delta_N)]))$. Such choice implies 
\bqn
\GW(\A^v_N(V[\mu^N(k\Delta_N)]),V[\mu^N(k\Delta_N)])\leq W(\A^v_N(V[\mu^N(k\Delta_N)]),V[\mu^N(k\Delta_N)])\leq |\mu^N(k\Delta_N)|\Delta^v_N.
\eqnl{e-Averror}

We now estimate recursively $\gw(\mu^N((k+1)\Delta_N),\nu^N((k+1)\Delta_N))$ starting from $\gw(\mu^N(k\Delta_N),\nu^N(k\Delta_N))$. Consider the PVFs 
$$V_1:=\A^v_N(V[\mu^N(k\Delta_N)]),\qquad  V_2=\A^v_N(V[\nu^N(k\Delta_N)]),$$ and the operator $\GW(V_1,V_2)$. By definition of $\GW$, there exist a choice $\tilde V_1\leq V_1,\tilde V_2\leq V_2$, and a transference plan $p\in  P(\tilde V_1,\tilde V_2)$ with the two following properties: one one side, it holds $\GW(V_1,V_2)=\int_{T\R^n\times T\R^n} |v-w|\, d p(x,v,y,w)$; on the other side, denoting with $\tilde\mu_i=\pi_1\# \tilde V_i$, it holds $\gw(\mu^N(k\Delta_N),\nu^N(k\Delta_N))=|\mu^N(k\Delta_N)-\tilde\mu_1|+W(\tilde\mu_1,\tilde\mu_2)+|\nu^N(k\Delta_N)-\tilde\mu_2|$, and $W(\tilde\mu_1,\tilde\mu_2)$ is realized by the transference plan $\pi_{13}\# p\in  P(\tilde\mu_1,\tilde\mu_2)$.

Consider now the following corresponding decomposition: write $\mu^N(k\Delta_N)=\sum_i m_i\delta_{x_i} +\sum_j n_j\delta_{y_j}$ and $\nu^N(k\Delta_N)=\sum_l p_l \delta_{z_i}+\sum_j n_j \delta_{t_j}$, where $\tilde\mu_1=\sum_j n_j\delta_{y_j}$, $\tilde\mu_2=\sum_j n_j \delta_{t_j}$, and the optimal transference plan $\pi_{13}\# p$ sends each $n_j\delta_{y_j}$ to each $n_j\delta_{t_j}$. Decompose accordingly the PFV as follows: 
$$V_1=\sum_{i,k}m_{ik}\delta_{(x_i,v_k)}+\tilde V^1\mbox{~~~with~}\tilde V^1:=\sum_{j,k}n_{jk}\delta_{(y_j,v_k)},$$
and similarly 
$$V_2=\sum_{l,k}p_{lk}\delta_{(z_i,v'_k)}+\tilde V^2\mbox{~~~with~}\tilde V^2:=\sum_{j,k}n_{jk}\delta_{(t_j,v'_k)},$$
with the additional requirement that the optimal transference plan $p\in P(\tilde V^1,\tilde V^2)$ sends $n_{jk}\delta_{(y_j,v_k)}$ to $n_{jk}\delta_{(t_j,v'_k)}$.

By definition of $\mu^N$, it then holds $\mu^N((k+1)\Delta_N)=\sum_{ik}m_{ik}\delta_{x_i+\Delta_N v_k}+\sum_{jk}n_{jk}\delta_{y_j+\Delta_N v_k}$, and similarly for $\nu^N((k+1)\Delta_N)$. Estimate the distance $\gw(\mu^N((k+1)\Delta_N),\nu^N((k+1)\Delta_N))$ by choosing the first component for mass removal and the second one for transport. It then holds
\bqn
\gw(\mu^N((k+1)\Delta_N),\nu^N((k+1)\Delta_N))&\leq& \sum_{ik}|m_{ik}|+\sum_{jk} n_{jk} |y_j+\Delta_Nv_k-t_j-\Delta_N v'_k|+ \sum_{lk}|p_{lk}|\leq\nn
&\leq& \sum_im_i +\sum_j n_j|y_j-t_j|+\Delta_N\sum_{jk}n_{jk}|v_k-v'_k|+\sum_l p_l=\nn
&=&\gw(\mu^N(k\Delta_N),\nu^N(k\Delta_N))+\Delta_N \GW(V_1,V_2).
\eqnl{e-gwk1}
We now need to compare $\GW(V_1,V_2)$ with $\GW(V[\mu^N(k\Delta_N)],V[\nu^N(k\Delta_N)])$. Denote with $\tilde V_1,\tilde V_2$ the decomposition and $p\in P(\tilde V_1,\tilde V_2)$ the transference plan in \r{e-GW} realizing $\GW(V[\mu^N(k\Delta_N)],V[\nu^N(k\Delta_N)])$. Observe that $\pi_1\#\tilde V_1\leq \pi_1\#(V[\mu^N(k\Delta_N)])=\mu^N(k\Delta_N)$ that is a finite sum of Dirac deltas, and the same holds for $\pi_1\#\tilde V_2$. Thus, $p$ can be decomposed as follows 
\bqn
p:=\sum_{i,k}p_{ik}(v,w) \delta_{x_i,y_k},
\eqnl{e-p}
where each $p_{ik}$ is a transference plan on $T_{x_i}\R^n\times T_{y_k}\R^n$.

We are now ready to define a decomposition $\bar V_1\leq V_1,\bar V_2 \leq V_2$ and a transference plan $q\in P(\bar V_1,\bar V_2)$ to estimate $\GW(V_1,V_2)$ from above. For each transference plan $p_{ik}$, define $$q_{ik}:=\sum_{jl}p_{ik}((v_j+Q')\times (w_l+Q'))\delta_{v_j,w_l},$$ where the $v_j,w_l$ are the equispaced discretized points on $T_x\R^n,T_y\R^n$, respectively, and $Q'$ is defined  in Definition \ref{d-discrete}. Then define
\bqn
q:=\sum_{i,k}q_{ik} \delta_{x_i,y_k}.
\eqnl{e-q}

Define now $\bar V_1:=\pi_{12}\#q$. By the definition of $\A^v_N$ in \r{e-defAv} and of $q$ in \r{e-q}, it is easy to prove that it holds
\bqn
\bar V_1=\pi_{12}\#q=\A^v_N(\pi_{12}\#p)=\A^v_N(\tilde V_1)\leq A^v_N(V[\mu^N(k\Delta_N)])=V_1.
\eqnn
One can equivalently prove that it holds $\bar V_2:=\pi_{34}\#q\leq V_2$. Moreover, it holds $\pi_{13}\#q=\pi_{13}\#p\in P^{opt}(\mu^N(k\Delta_N)],\nu^N(k\Delta_N))$, where we used that $p$ is a minimizer of $\GW(V[\mu^N(k\Delta_N)],V[\nu^N(k\Delta_N)])$. Then, the decomposition $\bar V_1,\bar V_2$ with the transference plan $q$ is admissible in the right hand side of \r{e-GW}. 
It then holds
\bqn
\GW(V_1,V_2)\leq \int_{T\R^n\times T\R^n}|v-w|\,dq(x,v,y,w)=\sum_{ij}\int_{T_{x_i}\R^n\times T_{y_k}\R^n}|v-w|\,dq_{ik}(v,w).
\eqnl{e-costoq}
We estimate each term by following the definition of $q_{ik}$, as follows
\bqn
\int_{T_{x_i}\R^n\times T_{y_k}\R^n}|v-w|\,dq_{ik}(v,w)&=&
\sum_{jl}|v_j-w_l|\,dq_{ik}(v,w)=\sum_{jl}|v_j-w_l| p_{ik}((v_j+Q')\times (w_l+Q'))\leq\nn
&&\sum_{jl}(|v-w|+2\mathrm{diam}(Q')) (p_{ik})_{|_{((v_j+Q')\times (w_l+Q'))}}.
\eqnn
Here, we used that $|v_j-w_k|\leq |v_j-v|+|v-w|+|v-w_k|$ for any $v\in v_j+Q'$ and $w\in w_l+Q'$. Observe now that, by decomposition, it holds 
\bqn
\sum_{jl}(|v-w|+2\mathrm{diam}(Q')) (p_{ik})_{|_{((v_j+Q')\times (w_l+Q'))}}=\int |v-w| \, dp_{ik}+2\mathrm{diam}(Q')|p_{ik}|.
\eqnn
Summing over $ik$ in \r{e-costoq} and recalling the definition of $p$ in \r{e-p}, we have
\bqn
\GW(V_1,V_2)\leq \int |v-w|\,dp+2\mathrm{diam}(Q') |p|\leq \GW(V[\mu^N(k\Delta_N)],V[\nu^N(k\Delta_N)])+2\sqrt{n}\Delta_N |\mu^N(k\Delta_N)|.
\eqnl{e-GWV1V2}
We used here that $\mathrm{diam}(Q')\leq \sqrt{n}\Delta_N$ as a consequence of Definition \ref{d-discrete}, as well as the fact that it holds $$|p|=|\pi_{12}\# p|\leq |V[\mu^N(k\Delta_N)]|=|\mu^N(k\Delta_N)|\leq C_6.$$

Going back to \r{e-gwk1}, and using the Lipschitz continuity hypothesis {\bf (V2)} in \r{e-GWV1V2}, it holds 
\bqn
\gw(\mu^N((k+1)\Delta_N),\nu^N((k+1)\Delta_N))\leq (1+K) \gw(\mu^N(k\Delta_N),\nu^N(k\Delta_N))+2\sqrt{n} (\Delta_N)^2 C_6
\eqnl{e-Prop2}

Clearly, this estimate implies $$\gw(\mu^N(t),\nu^N(t))\leq e^{Kt} \gw(\mu^N(0),\nu^N(0))+2\sqrt{n} (\Delta_N)^2 C_6\frac{e^{Kt}-1}{K},$$ as long as $t=k\Delta_N$. By uniform Lipschitz continuity given by \r{e-ul}, and Proposition \ref{p-errorediscr}, this implies
\bqn
\gw(\mu^N(t),\nu^N(t))\leq e^{Kt} \gw(\mu_0,\nu_0)+2C_6\Delta_N+2\sqrt{n} (\Delta_N)^2 C_6\frac{e^{Kt}-1}{K}.
\eqnl{e-key}

We now use a density argument to pass to the limit. Consider the countable set 
$$\mathcal{D}=\Pg{\mu_0\in \M \mbox{~~s.t.~~}\mu_0=\sum_{i}m_i\delta_{x_i}, 0<m_i\in\Q,x_i\in\Q^n}.$$
 Choose a $\mu_0\in\mathcal{D}$, define the corresponding sequence $\mu^N_0$ and choose a subsequence $N_k$ such that $\mu^{N_k}_0$ uniformly converges to a solution $\mu(t)$ to \r{e-MDES}. Then choose $\mu_1\in\mathcal{D}$, consider the subsequence $\mu_1^{N_k}$ and choose a converging subsequence $\mu_1^{N_{k_l}}$. Repeat this diagonal argument for the countable set of initial data in $\mathcal{D}$, and observe that, passing to the limit in \r{e-key} for $N\to\infty$, it holds
\bqn
\gw(\mu(t),\nu(t))\leq e^{Kt}\gw(\mu_0,\nu_0)
\eqnn
for all $\mu_0,\nu_0\in \mathcal{D}$. Observe now that $\mathcal{D}$ is dense in $\M$, thus the continuous semigroup $\mu_0\to S_t\mu_0=\mu(t)$ can be uniquely extended from $\mathcal{D}$ to $\M$. \hfill$\square$

\subsection{Uniqueness} \label{s-uniqueness}

We now prove Theorem \ref{t-uniqueness}, i.e. uniqueness of a solution to \r{e-MDES} when a Dirac germ $\gamma$ is fixed. We first need to define compatibility of a semigroup for the dynamics \r{e-MDES}, that is the following.
\bdeff\label{d-compat2} Fix a PVF $V$ satisfying {\bf (V1)}, a final time $T>0$ and a Dirac germ $\gamma$ as in Definition \ref{d-germ}. A semigroup for \r{e-MDES} is said to be compatible with $\gamma$ if one has the following property: for each $R,M>0$ there exists $C(R,M)$ such that the space $\M^D_{R,M}:=\Pg{\mu\in\M^D\mbox{~~s.t.~~}\supp(\mu)\in B(0,R), |\mu|\leq M}$ satisfies
\bqn
\mbox{for all }t\in[0,\inf_{\mu\in\M^D_{R,M}}\eps(\mu)] \mbox{ one has} \sup_{\mu\in\M^D_R} \gw(S_t\mu,\gamma_\mu(t))\leq C(R,M) t^2.
\eqnl{e-compat2}
\edeff
Observe that this definition coincides with Definition \ref{d-compatible}, where the dynamics \r{e-MDE} is replaced by \r{e-MDES} and the metric $W$ is replaced by $\gw$.

The proof of Theorem \ref{t-uniqueness} is then based on the following Lemma.
\bl\label{l-bressan} Let $S$ be a Lipschitz semigroup and $\mu:[0,T]\to\M$ a Lipschitz continuous curve. It then holds
\bqn
\gw(S_t\mu(0),\mu(t))\leq e^{Ct}\int_0^t \liminf_{h\to0^+}\gw(S_h\mu(s),\mu(s+h))\,ds,
\eqnn
where $C$ is the Lipschitz constant in Property 3b of Definition \ref{d-semigroup}.
\el
\bproof The original proof in Banach spaces can be found in \cite[Thm 2.9]{bressan}. Its adaptation to metric spaces can be found in \cite[Appendix A]{MDE}. 
\eproof

We are now ready to prove Theorem \ref{t-uniqueness}.

{\it Proof of Theorem \ref{t-uniqueness}}. Assume to have two semigroups $S^1_t,S^2_t$ for solutions to \r{e-MDES}, both compatible with a given germ $\gamma$. Fix an initial data $\mu_0\in\M$ and $T\geq0$. By Property 3a in the Definition \ref{d-semigroup}, there exists $R>0$ such that $\supp(S^1_t\mu_0)\cup\supp(S^2_t\mu_0)\subset B(0,R)$ for all $t\in[0,T]$. By applying Lemma \ref{l-bressan} with $\mu(t)=S^2_t$, it holds
\bqn
\gw(S^1_t\mu_0,S^2_t\mu_0)\leq e^{Ct}\int_0^t \liminf_{h\to0^+}\gw(S^1_hS^2_s\mu_0,S^2_{s+h}\mu_0)\,ds,
\eqnl{e-sem}

Fix $s>0$ and define $\nu=S^2_s\mu_0$. By density of $\M^D$ in $\M$ with respect to the topology induced by $\gw$, for each $\eps>0$ there exists $\bar\nu\in\M^D$ such that $\gw(\nu,\bar\nu)<\eps$. Property \r{e-compat2} applied to both $S^1_h$ and $S^2_h$ implies the existence of a constant $C_1(R,M)$ such that
$$\gw(S^1_h\bar\nu,\gamma_{\bar\nu}(h))\leq C_1(R,M) h^2,\qquad \gw(S^2_h\bar\nu,\gamma_{\bar\nu}(h))\leq C_1(R,M) h^2.$$
It then holds
\bqn
&&\gw(S^1_hS^2_s\mu_0,S^2_{s+h}\mu_0)\leq \gw(S^1_h\nu,S^1_h \bar\nu)+
\gw(S^1_h\bar\nu,\gamma_{\bar\nu}(h))+\gw(\gamma_{\bar\nu}(h),S^2_h\bar\nu)+
\gw(S^2_h\bar\nu,S^2_{h}\nu)\leq\nn
&&2(e^{Ch}\eps+C_1(R,M)h^2).
\eqnn
Since $s$ has been chosen arbitrarily and $\eps$ is arbitrarily small for the density property, it holds 
$$\liminf_{h\to 0^+}\frac{1}{h} \gw(S^1_hS^2_s\mu_0,S^2_{s+h}\mu_0)=0$$
for each $s>0$. This estimate applied to \r{e-sem} implies $S^1_t\mu_0=S^2_t\mu_0$ for any $t\geq 0$ and any initial data $\mu_0$, i.e. $S^1=S^2$.
 \hfill$\square$
 
\bibliographystyle{plain}
\bibliography{SourcePVF}
\end{document}